\documentclass[12pt]{amsart}

\begin{document}

\title{ Index transforms with the product  of the  associated Legendre functions}

\author{S. Yakubovich}
\address{Department of Mathematics, Fac. Sciences of University of Porto,Rua do Campo Alegre,  687; 4169-007 Porto (Portugal)}

\maketitle

\markboth{S. Yakubovich}{Index transforms with Legendre  functions}

\begin{abstract}   Index transforms with the product of the associated Legendre functions are introduced.      Mapping properties  are investigated  in the Lebesgue  spaces.  Inversion formulas  are proved.    The  results are applied to solve a  boundary value problem  in a wedge for  a  third  order partial differential equation.  
\end{abstract}


{\bf MS Classification}:   44A15, 33C10, 33C45, 44A05\bigskip

\section{Introduction and preliminary results}

Let $\mu \in \mathbb{C},\  f(x), g(\tau),\ x,\tau \in \mathbb{R}_+$ be complex-valued functions.  In the present paper we will investigate mapping properties of the following index transforms \cite{yak}

$$F(\tau)= \sqrt\pi \Gamma(1+i\tau-\mu)  \Gamma(1-i\tau-\mu) \int_0^\infty  P^\mu_{i\tau} (\sqrt{1+x}) P^\mu_{- i\tau} (\sqrt{1+x}) $$

$$\times {f(x)\over \sqrt{1+x}} dx,\eqno(1.1)$$

$$G(x)= \sqrt{ {\pi\over 1+x}}  \int_0^\infty \Gamma(1+i\tau-\mu)  \Gamma(1-i\tau-\mu)   P^\mu_{i\tau} (\sqrt{1+x}) P^\mu_{- i\tau} (\sqrt{1+x}) $$

$$\times g(\tau) d\tau. \eqno(1.2)$$
Here $i$ is the imaginary unit, $\Gamma(z)$ is the Euler gamma-function and   $P^\mu_\nu(z)$ is the associated Legendre function of the first kind  
(see \cite{erd}, Vol. I, \cite{vir}).  Our goal is to study them mapping properties, prove inversion theorems and apply to solve a  boundary value  problem for a  higher order PDE.   Denoting the kernel of (1.1), (1.2) by

 $$\Phi(x, \tau) =   \sqrt { {\pi \over 1+x}} \Gamma(1+i\tau-\mu)  \Gamma(1-i\tau-\mu) P^\mu_{i\tau} (\sqrt{1+x}) P^\mu_{- i\tau} (\sqrt{1+x}) ,\eqno(1.3)$$
we will find for further use its representation in terms of Fourier cosine transform \cite{tit} and deduce an ordinary differential equation with polynomial coefficients, whose solution is $\Phi(x,\tau)$, employing  the so-called method of the Mellin-Barnes integrals, which is already being successfully applied by the author for other index transforms.  In fact, appealing  to \cite{prud}, Vol. III,  entry 8.4.41.50,  we find the following Mellin-Barnes integral representation for the kernel (1.3), namely,

$$ \Phi(x,\tau ) = {1\over 2\pi i} \int_{\gamma-i\infty}^{\gamma +i\infty} \frac {\Gamma(1-s+  i\tau)\Gamma(1-s -i\tau) \Gamma(1/2-s)\Gamma(s-\mu) }{\Gamma(1-s) \Gamma (1-s-\mu) } x^{-s} ds, \ x >0,\eqno(1.4)$$
where $\gamma$ is taken from the interval $( {\rm Re} \mu,\ 1/2)$. The absolute convergence of the integral (1.4)   follows immediately from  the Stirling asymptotic formula for the gamma- function \cite{erd}, Vol. I, because for all $\tau \in \mathbb{R}$
$$\frac {\Gamma(1-s+  i\tau)\Gamma(1-s -i\tau) \Gamma(1/2-s)\Gamma(s-\mu) }{\Gamma(1-s) \Gamma (1-s-\mu) } $$

$$  =  O \left(e^{-\pi |s|} |s|^ {-1/2}\right),\  |s| \to \infty.\eqno(1.5)$$
Moreover, it can be differentiated under the integral sign any number of times due to the absolute and uniform convergence by $x \ge x_0 >0$.   Our method of investigation of the index transforms (1.1), (1.2) is based on   the Mellin transform technique developed in \cite{yal}.   Precisely, the Mellin transform is defined, for instance, in  $L_{\nu, p}(\mathbb{R}_+),\ 1 \le  p \le 2$ (see details in \cite{tit}) by the integral  
$$f^*(s)= \int_0^\infty f(x) x^{s-1} dx,\eqno(1.6)$$
 being convergent  in mean with respect to the norm in $L_q(\nu- i\infty, \nu + i\infty),\ \nu \in \mathbb{R}, \   q=p/(p-1)$.   Moreover, the  Parseval equality holds for $f \in L_{\nu, p}(\mathbb{R}_+),\  g \in L_{1-\nu, q}(\mathbb{R}_+)$
$$\int_0^\infty f(x) g(x) dx= {1\over 2\pi i} \int_{\nu- i\infty}^{\nu+i\infty} f^*(s) g^*(1-s) ds.\eqno(1.7)$$
The inverse Mellin transform is given accordingly
 $$f(x)= {1\over 2\pi i}  \int_{\nu- i\infty}^{\nu+i\infty} f^*(s)  x^{-s} ds,\eqno(1.8)$$
where the integral converges in mean with respect to the norm  in   $L_{\nu, p}(\mathbb{R}_+)$
$$||f||_{\nu,p} = \left( \int_0^\infty  |f(x)|^p x^{\nu p-1} dx\right)^{1/p}.\eqno(1.9)$$
In particular, letting $\nu= 1/p$ we get the usual space $L_p(\mathbb{R}_+)$.   

We begin with

{\bf Lemma 1}. {\it Let $x,  \tau \in \mathbb{R}_+,\  {\rm Re}\mu < 1/2$. Then the kernel $(1.3)$ has the following integral representation in terms of  Fourier cosine transform of the associated Legendre function} 
$$ \Phi(x,\tau)  =     \Gamma\left( {3\over 2}-\mu\right) e^{\pi i\mu} \int_0^\infty  {\cos(\tau u) \over \cosh^{1/2} (u/2)  (x+  \cosh^2(u/2) )^{3/4} } $$

$$\times  P_{1/2}^\mu \left(  {  x +2 \cosh^2(u/2) \over 2\cosh(u)   (  x + \cosh^2 (u/2) )^{1/2} } \right) du.\eqno(1.10) $$

\begin{proof}     In fact, appealing  to the reciprocal formulae via the Fourier cosine transform (cf. formula (1.104) in \cite{yak}) 
$$\int_0^\infty  \Gamma\left(1-s  + i\tau\right)  \Gamma\left(1-s  - i\tau \right)  \cos( \tau y) d\tau
= {\pi\over 2^{2(1-s)}}  {\Gamma(2(1-s)) \over \cosh^{2(1-s)}(y/2)},\ {\rm Re}\ s  < 1,\eqno(1.11)$$
$$  \Gamma\left(1-s + i\tau \right)  \Gamma\left(1-s  - i\tau\right)  
=   { \Gamma(2(1-s))  \over 2^{1-2s}}  \int_0^\infty   {\cos(\tau y)  \over \cosh^{2(1-s)} (y/2)} \ dy,\eqno(1.12)$$ 
we replace  the gamma-product $ \Gamma\left(1-s + i\tau \right)  \Gamma\left(1-s - i\tau\right)$ in the integral (1.4)  by its integral representation (1.12) and change the order of integration via Fubini's theorem.  Then, employing the duplication formula for the gamma-function \cite{erd}, Vol. I,   we derive 
$$ \Phi (x, \tau)  = {1\over 2 \pi\sqrt \pi  i} \int_0^\infty   {\cos(\tau u) \over \cosh^2 (u/2)}  $$

$$\times \int_{\gamma-i\infty}^{\gamma +i\infty} \frac { \Gamma(s-\mu) \Gamma(1/2- s)\Gamma(3/2-s) }{ \Gamma (1-s-\mu) } \left( {x\over  \cosh^2(u/2)} \right)^{-s}  ds du .\eqno(1.13) $$
In the meantime,  the  inner integral with respect to $s$ can be calculated, employing the Parseval equality (1.7), properties of the Mellin transform (1.6) and relations (8.4.3.1), (8.4.46.7) in \cite{prud}, Vol. III.   Thus we obtain

$${1\over 2\pi i} \int_{\gamma-i\infty}^{\gamma +i\infty} \frac { \Gamma(s-\mu) \Gamma(1/2- s)\Gamma(3/2-s) }{ \Gamma (1-s-\mu) } \left( {x \over \cosh^2(u/2)} \right)^{-s}  ds$$

$$  =  \left( {x \over  \cosh^2(u/2) } \right)^{-\mu} \int_0^\infty y^{1/2 -\mu} e^{- (x  \cosh^{-2}(u/2) +1) y} \ \Psi\left( {1\over 2}-\mu, 2; y\right) dy,\eqno(1.14)$$
where $\Psi(a, b; z)$ is the Tricomi function \cite{erd}, Vol. I.  But entry 7.11.4.9 in \cite{prud}, Vol. III permits to express Tricomi's function in (1.14) in terms of the Bateman function $k_\nu(z)$ \cite{erd}, Vol. II, and we have

$$ \left( {x \over  \cosh^2(u/2)}  \right)^{-\mu} \int_0^\infty y^{1/2 -\mu} e^{- (x  \cosh^{-2}(u/2) +1) y} \ \Psi\left( {1\over 2}-\mu, 2; y\right) dy$$

$$= \Gamma(3/2+\mu)  \left( {x \over  \cosh^2(u/2) } \right)^{-\mu} \int_0^\infty y^{-1/2 -\mu} e^{- (x  \cosh^{-2}(u/2) +1/2) y} \  k_{1+2\mu} \left({y\over 2} \right) dy.$$
Meanwhile, the latter integral is calculated in \cite{prud}, Vol. III, entry 2.14.2.4 in terms of the Gauss hypergeometric function \cite{erd}, Vol. I and it gives after slight simplifications 

$$ \Gamma(3/2+\mu)  \left( {x\over   \cosh^2(u/2)}  \right)^{-\mu} \int_0^\infty y^{-1/2 -\mu} e^{- (x  \cosh^{-2}(u/2) +1/2) y} \  k_{1+2\mu} \left({y\over 2} \right) dy$$

$$=  4^{\mu} \sqrt \pi  \left( {x \over  \cosh^2\left(u/2 \right) } \right)^{-\mu} \frac{ \Gamma (3/2-\mu) }{ \Gamma(1-\mu)}
 \ {}_2F_1 \left( {1\over 2}-\mu, \  {3\over 2}-\mu ;\  1- 2\mu; \ -  {x \over  \cosh^2\left(u/2 \right) }\right).$$
Finally, appealing to relation (7.3.1.70) in \cite{prud}, Vol. III, we express the value of the Gauss hypergeometric function in terms of the associated Legendre function.  Precisely, we find 

$$4^{\mu} \sqrt \pi  \left( {x \over  \cosh^2\left(u/2 \right) } \right)^{-\mu} \frac{ \Gamma (3/2-\mu) }{ \Gamma(1-\mu)}
 \ {}_2F_1 \left( {1\over 2}-\mu, \  {3\over 2}-\mu ;\  1- 2\mu; \ -  {x  \over \cosh^2\left(u/2 \right) }\right)$$

$$= \sqrt\pi \Gamma\left( {3\over 2}-\mu\right) e^{\pi i\mu} \left( 1+  x  \cosh^{-2}\left({u\over 2} \right)\right)^{-3/4} 
P_{1/2}^\mu \left(  { x+ 2\cosh ^2(u)  \over 2 \cosh(u)  (  x +  \cosh^{2} (u/2) )^{1/2} } \right).$$
Hence substituting this value in the right-hand side of (1.14) and recalling (1.13),  we arrive at the formula (1.10), completing the proof of Lemma 1. 

\end{proof}

{\bf Lemma 2}.    {\it The function  $ \Phi(x, \tau),\ (x,\tau ) \in \mathbb{R}_+ \times  \mathbb{R}_+ $ given by formula $(1.3)$ satisfies  the following third order differential equation with polynomial   coefficients}
$$2x^2(1+x)  {d^3 \Phi\over dx^3} + x(11x+ 6) {d^2 \Phi \over dx^2}  $$

$$+   \left( 2(1-\mu^2) + x( 11 + 2\tau^2)  \right)  {d\Phi \over dx}  + \left( 1 + \tau^2\right) \Phi = 0. \eqno(1.15)$$

\begin{proof}  As it was mentioned above, the asymptotic behavior (1.5) of the integrand in (1.4) permits a differentiation under the integral sign any number of times.  Hence employing the reduction formula for the gamma- function \cite{erd}, Vol. I , we derive 

$$\left( x {d\over dx} \right)^2  \Phi(x, \tau) =   {1\over 2\pi i} \int_{\gamma-i\infty}^{\gamma +i\infty} \frac { s^2 \Gamma(1-s+  i\tau)\Gamma(1-s -i\tau) \Gamma(1/2-s)\Gamma(s-\mu) }{\Gamma(1-s) \Gamma (1-s-\mu) }  x^{-s} ds$$

$$=   {1\over 2\pi i} \int_{\gamma-i\infty}^{\gamma +i\infty} \frac {  \Gamma(2-s+  i\tau)\Gamma(2-s -i\tau) \Gamma(1/2-s)\Gamma(s-\mu) }{\Gamma(1-s) \Gamma (1-s-\mu) }  x^{-s} ds$$

$$ -  \left(1+ \tau^2\right)  \Phi(x, \tau) - 2 x {d\Phi\over dx}.\eqno(1.16)$$ 
Meanwhile, with a simple change of variable

$$ {1\over 2\pi i} \int_{\gamma-i\infty}^{\gamma +i\infty} \frac {  \Gamma(2-s+  i\tau)\Gamma(2-s -i\tau) \Gamma(1/2-s)\Gamma(s-\mu) }{\Gamma(1-s) \Gamma (1-s-\mu) }  x^{-s} ds$$

$$=  {1\over 2\pi i} \int_{\gamma-1-i\infty}^{\gamma-1 +i\infty} \frac {  \Gamma(1-s+  i\tau)\Gamma(1-s -i\tau) \Gamma(-1/2-s)\Gamma(s+1-\mu) }{\Gamma(-s) \Gamma (-s-\mu) }  x^{-s-1} ds$$

$$= -  {1\over 2\pi i} \int_{\gamma-1-i\infty}^{\gamma-1 +i\infty} \frac {s  (s^2-\mu^2) \Gamma(1-s+  i\tau)\Gamma(1-s -i\tau) \Gamma(1/2-s)\Gamma(s-\mu) }{(1/2+s) \Gamma(1-s) \Gamma (1-s-\mu) }  x^{-s-1} ds$$

$$= -  {1\over 2\pi i} \int_{\gamma-1-i\infty}^{\gamma-1 +i\infty} \frac { (s^2-\mu^2) \Gamma(1-s+  i\tau)\Gamma(1-s -i\tau) \Gamma(1/2-s)\Gamma(s-\mu) }{ \Gamma(1-s) \Gamma (1-s-\mu) }  x^{-s-1} ds$$

$$+  {1\over 4\pi i} \int_{\gamma-1-i\infty}^{\gamma-1 +i\infty} \frac {  (s^2-\mu^2) \Gamma(1-s+  i\tau)\Gamma(1-s -i\tau) \Gamma(1/2-s)\Gamma(s-\mu) }{(1/2+s) \Gamma(1-s) \Gamma (1-s-\mu) }  x^{-s-1} ds.$$
Therefore, moving the contour of  two latter integrals to the right via Cauchy's theorem, multiplying  by $\sqrt x$, differentiating again and using (1.16), we obtain

$${d\over dx} \left[ \sqrt x \left[ \left( x {d\over dx} \right)^2  \Phi  +  \left(1+ \tau^2\right)  \Phi + 2 x {d\Phi \over dx} \right.\right.$$

$$\left.\left.  +{1\over x}  \left( x {d\over dx} \right)^2  \Phi - {\mu^2\over x}  \Phi\right] \right]  =  - {x^{-3/2} \over 2} \left[  \left( x {d\over dx} \right)^2  \Phi  -\mu^2  \Phi \right].$$
Hence we arrive at the following operator equation

$$2  \left( x {d\over dx} \right)^3  \Phi + 5 \left( x {d\over dx} \right)^2  \Phi  + 2(2+\tau^2) \left( x {d\over dx} \right)  \Phi  + 2 {d\over dx} \left( \left( x {d\over dx} \right)^2 -\mu^2 \right)  \Phi + (1+\tau^2) \Phi =0.$$ 
Thus fulfilling the differentiation, we end up with (1.15).    Lemma 2 is proved. 

\end{proof}

\section {Boundedness  and inversion properties for the index transform (1.1)}

In this section we will investigate the mapping properties of the index transform (1.1).  In fact, denoting by $C_b(\mathbb{R}_+)$ the space of bounded continuous functions, we establish

{\bf Theorem 1.}   {\it Let ${\rm Re}\mu < 1/2$. The index transform  $(1.1)$  is well-defined as a  bounded operator $F: L_{1-\nu, 1} \left(\mathbb{R}_+\right) \to C_b (\mathbb{R}_+),\  \nu \in ( {\rm Re}\mu,\  1/2)$ and the following norm inequality takes place
$$||F  f||_{C_b(\mathbb{R}_+)} \equiv \sup_{\tau \in \mathbb{R}_+} | (F f)(\tau)| \le C_{\mu,\nu}  ||f||_{1-\nu,1},\eqno(2.1)$$
where

$$C_{\mu,\nu} =  {2^{-2\nu } \over  \pi\sqrt \pi } B\left(1-\nu,\  1-\nu  \right) \int_{\nu-i\infty}^{\nu +i\infty} \left| \frac {  \Gamma(3/2-s)  \Gamma(1/2-s)\Gamma(s-\mu) }{\Gamma (1-s-\mu) }  ds\right| \eqno(2.2)$$
and $B(a,b)$ is the Euler beta-function \cite{erd}, Vol. 1.  Moreover,  $(F f)(\tau) \to 0,\  \tau \to \infty.$  Besides, if, in addition,  $ f \in   L_{1-\nu, p}(\mathbb{R}_+),\  \varphi (x) \in L_{1-\nu,p}(\mathbb{R}_+),\ 1 < p\le 2$,  $ \nu \in \left( \max(0, {\rm Re} \mu),\  1/(2q) \right),\  {\rm Re} \mu<  1/(2q),\ q= p/(p-1)$,  where

$$\varphi(x)=  {1\over 2\pi i} \int_{\nu -i\infty}^{\nu  +i\infty} \frac { \Gamma(s -\mu)\Gamma(s) \Gamma(1/2-s)}{\Gamma(s-1/2) \Gamma(1 -s) \Gamma (1-s-\mu) } f^*(1 -s) x^{-s} ds,\eqno(2.3)$$
then  
$$(F f)(\tau) =   {\sqrt{\pi} \over\cosh(\pi\tau)} \int_0^\infty   K_{i\tau} \left(\sqrt x\right) \left[  I_{i\tau} \left(\sqrt x\right) +  I_{- i\tau} \left( \sqrt x\right)\right] \varphi(x)\  dx\eqno(2.4)$$
and  integrals $(2.3),\ (2.4)$ converge absolutely.}

\begin{proof}  Indeed, recalling (1.4), (1.9), (1.12), using the duplication formula for the gamma-function  and calculating an elementary integral with the hyperbolic function, we derive

$$ \left|(F f)(\tau) \right| \le \int_0^\infty \left| \Phi (x,\tau) \right|  |f(x)| dx $$

$$\le {1 \over \pi\sqrt\pi } \int_{\nu-i\infty}^{\nu  +i\infty}  \left| \frac {  \Gamma(3/2-s)  \Gamma(1/2-s)\Gamma(s-\mu) }{\Gamma (1-s-\mu) } ds\right|  $$

$$\times \int_0^\infty  |f(x)|  x^{-\nu} dx  \int_0^\infty   {dy \over \cosh^{2(1-\nu)} (y)} $$

$$=  { 2^{-2\nu} \over \pi\sqrt\pi } \  B\left(1-\nu, 1-\nu\right) \int_{\nu-i\infty}^{\nu  +i\infty} \left| \frac {  \Gamma(3/2-s)  \Gamma(1/2-s)\Gamma(s-\mu) }{\Gamma (1-s-\mu) }  ds\right|     ||f||_{1-\nu,1}.$$
The latter estimate proves (2.1).   Furthermore, Fubini's theorem and the definition of the Mellin transform (1.6) yield the integral representation of $(Ff)(\tau)$ in terms of the Fourier cosine transform of integrable function, namely,

$$ (F f)(\tau) =    {1 \over \pi i \sqrt\pi }  \int_0^\infty   {\cos(2 \tau y)  \over \cosh^{2} (y)}  \int_{\nu-i\infty}^{\nu  +i\infty}  \frac {  \Gamma(3/2-s)  \Gamma(1/2-s)\Gamma(s-\mu) }{\Gamma (1-s-\mu) }$$

$$\times   f^*(1-s) \cosh^{2s} (y) ds dy.$$ 
Hence it tends to zero, when $\tau \to \infty$ via the Riemann-Lebesgue lemma.  Moreover, the Parseval equality (1.7) and (1.4) give the representation 

$$ (Ff)(\tau) = {1\over 2\pi i} \int_{\nu-i\infty}^{\nu  +i\infty} \frac {\Gamma(1-s+  i\tau)\Gamma(1-s -i\tau) \Gamma(1/2-s)\Gamma(s-\mu) }{\Gamma(1-s) \Gamma (1-s-\mu) }$$

$$\times   f^*(1-s) ds.\eqno(2.5)$$
Hence, employing  the Mellin-Barnes representation for the product of modified Bessel functions of the third kind (cf. relation (8.4.23.23) in \cite{prud}, Vol. III)

$${\sqrt{\pi} \over \cosh(\pi\tau)} K_{i\tau} \left( \sqrt x\right) \left[  I_{i\tau} \left( \sqrt x\right) +  I_{- i\tau} \left( \sqrt x \right)\right] $$

$$={1\over 2\pi i} \int_{\nu -i\infty}^{\nu  +i\infty} \frac { \Gamma( 1/2-s) }{\Gamma(1-s)} \   \Gamma(s+  i\tau)\Gamma(s -i\tau) x^{-s} ds,\ x >0$$
 and using again the Parseval equality (1.7), we find that (2.5) becomes the Lebedev index transform with the product of the modified Bessel functions \cite{square} given by formula (2.4),  where $\varphi(x)$ is  defined by integral (2.3). Its  absolute convergence for each $x >0$  can be verified via  H\"older's  inequality 

$$ \int_{\nu-i\infty}^{\nu +i\infty} \left| \frac { \Gamma(s -\mu)\Gamma(s) \Gamma(1/2-s)}{\Gamma(s-1/2) \Gamma(1 -s) \Gamma (1-s-\mu) } f^*(1 -s) x^{-s} ds\right| $$

$$\le x^{-\nu} \left( \int_{\nu -i\infty}^{\nu  +i\infty} \left|f^*(1-s)\right|^q |ds|\right)^{1/q}$$

$$\times \left( \int_{\nu -i\infty}^{\nu  +i\infty} \left|\frac { \Gamma(s -\mu)\Gamma(s) \Gamma(1/2 -s)}{\Gamma(s-1/2) \Gamma(1-s) \Gamma (1-s-\mu) } \right|^p |ds|\right)^{1/p}  < \infty,\ q= {p\over p-1}.$$
The convergence of the latter integral by $s$ is justified,  recalling  the Stirling formula for the asymptotic behavior of the gamma-function, which gives

$$ \frac { \Gamma(s-\mu)\Gamma(s) \Gamma(1/2-s )}{\Gamma(s-1/2) \Gamma(1-s) \Gamma (1-s-\mu) } 
= O\left( |s|^{2\nu-1}\right),\ |s| \to \infty,$$
and $\nu$ is chosen from the interval  $\left(\max (0, {\rm Re} \mu),\ 1/(2q) \right),\   {\rm Re} \mu <  1/(2q)$.    In order to establish the absolute convergence of the integral (2.4), we use the assumption $\varphi (x) \in L_{1-\nu,p}(\mathbb{R}_+)$ and invoke the asymptotic formulae for the modified Bessel functions \cite{erd}, Vol. 2 for fixed $\tau \in \mathbb{R}$, namely, 

$$K_{i\tau} \left(\sqrt {x}\right) \left[  I_{i\tau} \left(\sqrt { x}\right) +  I_{- i\tau} \left(\sqrt { x}\right)\right] = O(\log x),\  x\to 0+, $$

$$ K_{i\tau} \left(\sqrt {x}\right) \left[  I_{i\tau} \left(\sqrt { x}\right) +  I_{- i\tau} \left(\sqrt { x}\right)\right] = O\left( x^{-1/2}\right),\ x \to \infty.$$
Then accordingly,   

$$\int_0^\infty \left| K_{i\tau} \left(\sqrt x\right) \left[  I_{i\tau} \left(\sqrt x\right) +  I_{- i\tau} \left( \sqrt x\right)\right] \varphi(x)\right|   dx$$

$$\le ||\varphi||_{1-\nu,p} \left(\int_0^\infty \left| K_{i\tau} \left(\sqrt x\right) \left[  I_{i\tau} \left(\sqrt x\right) +  I_{- i\tau} \left( \sqrt x\right)\right] \right|^q x^{\nu q-1}   dx\right)^{1/q} $$

$$= ||\varphi||_{1-\nu,p} \left(\int_0^1  O\left ( x^{\nu q-1}\  \log^q x \right)  dx  +  \int_1^\infty  O\left ( x^{q(\nu-1/2)-1} \right)  dx\right)^{1/q} < \infty. $$
Theorem 1 is proved. 

\end{proof}

Writing  (2.4) in the form

$$(F f)(\tau) =   {2 \sqrt{\pi} \over\cosh(\pi\tau)} \int_0^\infty   K_{i\tau} \left(x\right) \left[  I_{i\tau} \left(x\right) +  I_{- i\tau} \left( x\right)\right] \varphi(x^2) x  dx,\eqno(2.6)$$
we  appeal to the Lebedev expansion theorem in \cite{square},  which implies the following representation of the antiderivative 

$$ \int_{x}^{\infty}  \varphi\left(y^2\right)  y dy  = {1\over \pi^2\sqrt \pi} \int_0^\infty \tau \sinh(2\pi\tau) K_{i\tau}^2(x) (F f)(\tau) d\tau,\ x >0,\eqno(2.7)$$
which holds under condition $x \varphi(x^2) \in L_1\left((0,1); \ x^{-1/2} dx\right) \cap  L_1\left((1,\infty); \  x^{1/2} dx\right)$. By straightforward substitutions we see that this condition is equivalent to  (cf. (1.9)) $\varphi \in L_{5/4,1}(1,\infty) \cap L_{3/4,1}(0,1)$.   Hence (2.7) yields the equality  

$$ \int_{x}^{\infty}  \varphi\left(y\right)  dy  = {2\over \pi^2\sqrt \pi} \int_0^\infty \tau \sinh(2\pi\tau) K_{i\tau}^2(\sqrt x) (F f)(\tau) d\tau.\eqno(2.8)$$
Then, appealing to relation (8.4.23.27) in \cite{prud}, Vol. III and  differentiating both  sides of (2.8)  with respect to $x$, we find

$$\varphi(x) = - {1 \over 2\pi^3 i} {d\over dx}  \int_0^\infty \tau \sinh(2\pi\tau)  (F f)(\tau) $$

$$\times \int_{\nu -i\infty}^{\nu  +i\infty} \frac { \Gamma(s)\Gamma(s+ i\tau) \Gamma(s-i\tau)}{\Gamma(1/2 +s)} x^{-s} ds d\tau,\eqno(2.9)$$
where \  $\nu > 0.$  Our goal now is to motivate the differentiation under integral sign in the right-hand side of (2.9).   To do this, we first appeal  to the Parseval equality (1.7) and relation (8.4.23.3) in \cite{prud}, Vol. III in order to rewrite the integral with respect to $s$ in (2.10) as follows

$${1\over 2\pi i} \int_{\nu -i\infty}^{\nu  +i\infty} \frac { \Gamma(s)\Gamma(s+ i\tau) \Gamma(s-i\tau)}{\Gamma(1/2 +s)} x^{-s} ds ={1\over \sqrt \pi}  
\int_0^\infty e^{-y/2- x/y} K_{i\tau} \left({y\over 2}\right) {dy\over y} .$$
Hence after substitution of the right-hand side of the latter equality into (2.9) a formal differentiation means

$$\varphi(x) = {1 \over \pi^2 \sqrt \pi}  \int_0^\infty \tau \sinh(2\pi\tau)  (F f)(\tau)    \int_0^\infty e^{-y/2- x/y} K_{i\tau} \left({y\over 2}\right) {dy d\tau \over y^2}.\eqno(2.10)$$
It is allowed via the absolute and uniform convergence by $x \ge x_0 >0$ of the iterated integral (2.10), which can be verified under the integrability condition $Ff \in L_1(\mathbb{R}_+; \tau e^{3\pi \tau/2} d\tau)$ and the Lebedev inequality for the modified Bessel function \cite{yal}, p. 99

$$|K_{i\tau}(y)| \le {y^{-1/4} \over \sqrt{\sinh(\pi\tau)}},\ y, \tau  >0.\eqno(2.11)$$
Indeed, with the use of relation (2.3.16.1) in \cite{prud}, Vol. I we have 

$$\int_0^\infty \tau \sinh(2\pi\tau) \left| (F f)(\tau) \right|   \int_0^\infty e^{-y/2- x/y} \left| K_{i\tau} \left({y\over 2}\right) \right|  {dy d\tau \over y^2} $$

$$\le   \int_0^\infty \tau  e^{3\pi \tau/2} \left| (F f)(\tau) \right|   \int_0^\infty e^{-y - x_0/(2y)} {dy d\tau \over y^{9/4}} $$

$$=  2^{13/8}\  x_0^{-5/8} K_{5/4} \left( \sqrt{2x_0}\right) \int_0^\infty \tau  e^{3\pi \tau/2} \left| (F f)(\tau) \right| d\tau < \infty. $$
Moreover, one can change the order of integration in (2.10) due to Fubini's theorem to get

$$\varphi(x) = {1 \over \pi^2 \sqrt \pi}   \int_0^\infty e^{-y/2- x/y}   \int_0^\infty \tau \sinh(2\pi\tau)  K_{i\tau} \left({y\over 2}\right)  (F f)(\tau)   {d\tau dy\over y^2}.\eqno(2.12)$$
This equality guarantees the condition  $\varphi \in L_{\nu+1,1}\left(\mathbb{R}_+\right),\    \nu >  1/4.$  In fact, we have with (2.11)

$$||\varphi||_{1+\nu,1} = \int_0^\infty |\varphi(x)| x^\nu dx \le  {\Gamma(\nu+3/4) \Gamma(\nu-1/4)\  2^{\nu+3/4}   \over \pi^2 \sqrt \pi} $$

$$\times   \int_0^\infty \tau \sqrt{\sinh(\pi\tau)}\cosh(\pi\tau) \left|  (F f)(\tau) \right| d\tau < \infty. $$
Therefore,  recalling the Parseval equality (1.7), equality (2.12) can be written in the  form

$$x \varphi(x) = {1 \over 2\pi^3 i}  \int_{\nu -i\infty}^{\nu  +i\infty} \frac { \Gamma(1+s) } {\Gamma(1/2 +s)} x^{-s}  \int_0^\infty \tau \sinh(2\pi\tau)  (F f)(\tau)  $$

$$\times \Gamma(s+ i\tau) \Gamma(s-i\tau) d\tau ds.\eqno(2.13)$$
Besides,  applying  the Mellin transform to both sides (2.13) (cf. (2.3)), we  obtain after a slight simplification 

$$ \frac { \Gamma(s +1-\mu) \Gamma(-1/2-s)}{ \Gamma( -s) \Gamma (-s-\mu) } f^*( -s) = {1 \over \pi^2}  \int_0^\infty \tau \sinh(2\pi\tau)  (F f)(\tau)  $$

$$\times \Gamma(s+ i\tau) \Gamma(s-i\tau) d\tau.\eqno(2.14)$$
Hence, employing the reduction formula for the gamma-function and the inverse Mellin transform (1.8) under the integrability condition $s f^*(-s) \in L_1\left(\nu-i\infty,\ \nu+ i\infty\right)$, we deduce from (2.14)

$$- {1 \over 2\pi i}  \int_{\nu -i\infty}^{\nu  +i\infty} s f^*(-s) x^{-s} ds =  {1 \over 2\pi^3 i}  \int_{\nu -i\infty}^{\nu  +i\infty} \frac{ \Gamma(1 -s) \Gamma (-s-\mu) }  { \Gamma(s +1-\mu) \Gamma(-1/2-s)} x^{-s} $$

$$\times \int_0^\infty \tau \sinh(2\pi\tau)  (F f)(\tau)  \Gamma(s+ i\tau) \Gamma(s-i\tau) d\tau ds.\eqno(2.15)$$
The left-hand side of (2.15) can be treated with the use of the operational properties for the Mellin transform under corresponding conditions, and in the right-hand side we interchange the order of integration, justifying this passage below.  Consequently, we obtain for almost all $x >0$ the equality

$$ x {d\over dx} \left[  f\left({1\over x}\right) \right] = {1\over \pi^2}  \int_0^\infty \tau \sinh(2\pi\tau) S(x,\tau)  (F f)(\tau) d\tau,\eqno(2.16)$$ 
where the kernel $S(x,\tau)$ is defined in terms of the Mellin-Barnes integral

$$S(x,\tau)= {1 \over 2\pi i}  \int_{\nu -i\infty}^{\nu  +i\infty} \frac {\Gamma(s+ i\tau) \Gamma(s-i\tau)  \Gamma(1 -s) \Gamma (-s-\mu) }  { \Gamma(s +1-\mu) \Gamma(-1/2-s)} x^{-s} ds,\eqno(2.17)$$
and the vertical line $(\nu-i\infty,\ \nu+i\infty)$ in the complex plane $s$ separates left- and right-hand side simple poles of the gamma-functions under the assumption $0 < \nu <  \min( - {\rm Re} \mu, 1),\ \mu \notin \mathbb{Z}$.  

In the meantime, to justify the interchange of the order of integration in (2.15), we appeal again to the Parseval equality (1.7) and relation  (8.4.23.3) in \cite{prud}, Vol. III to write  the right-hand side of (2.15) in the form

$${1 \over 2\pi^3 i}  \int_{\nu -i\infty}^{\nu  +i\infty} \frac{ \Gamma(1 -s) \Gamma (-s-\mu) }  { \Gamma(s +1-\mu) \Gamma(-1/2-s)} x^{-s} $$

$$\times \int_0^\infty \tau \sinh(2\pi\tau)  (F f)(\tau)  \Gamma(s+ i\tau) \Gamma(s-i\tau) d\tau ds $$

$$= {1\over \pi^{5/2}} \int_0^\infty  e^{-y/2} h\left({x\over y} \right)   \int_0^\infty \tau \sinh(2\pi\tau) \ K_{i\tau} \left({y\over 2} \right)  (F f)(\tau)   {d\tau dy\over y},\eqno(2.18)$$
where 

$$h(x) =  {1 \over 2\pi i}  \int_{\nu -i\infty}^{\nu  +i\infty} \frac{ \Gamma(s+1/2) \Gamma(1 -s) \Gamma (-s-\mu) }  { \Gamma(s +1-\mu) \Gamma(-1/2-s)} x^{-s} ds\eqno(2.19) $$
and  $1/4 < \nu <  \min( - {\rm Re} \mu, 1)$. But for $\mu \notin \mathbb{Z}$ the integral (2.19) can be calculated in terms of the generalized hypergeometric functions by means of the Slater theorem \cite{prud}, Vol. III. Precisely, we find the value

$$h(x)=   { 3 \ \Gamma (-1- \mu)  \over  8 x \ \Gamma(2- \mu)} \ {}_2F_2 \left( {3\over 2},\ {5\over 2} ;\  2+\mu,\ 2-\mu;\ - {1\over x} \right)$$

$$+ \left( {x\over 4}\right)^\mu  {  \sqrt \pi\ \Gamma (1+ \mu)   \over   \Gamma(1-\mu) \Gamma( \mu-1/2)} \ {}_2F_2 \left( {1\over 2} - \mu,\ {3\over 2} - \mu ; \ - \mu,\ 1- 2\mu;\ - {1\over x} \right).\eqno(2.20)$$
Meanwhile,  the interchange of the order of integration in the right-hand side of (2.18) by Fubini's theorem is motivated by the following estimate (see (2.11), (2.19))

$$\int_0^\infty  e^{-y/2} \left| h\left({x\over y} \right) \right|   \int_0^\infty \tau \sinh(2\pi\tau) \  \left| K_{i\tau} \left({y\over 2} \right)  (F f)(\tau) \right|  {d\tau dy\over y} $$

$$\le  {x^{-\nu}  2^{1/4} \over \pi} \int_0^\infty  e^{-y/2} y^{\nu-5/4} dy  \int_{\nu -i\infty}^{\nu  +i\infty}  \left| \frac{ \Gamma(s+1/2) \Gamma(1 -s) \Gamma (-s-\mu) }  { \Gamma(s +1-\mu) \Gamma(-1/2-s)} ds \right |$$

$$\times  \int_0^\infty \tau \sqrt {\sinh(\pi\tau)} \cosh(\pi\tau)  \  \left|   (F f)(\tau) \right|  d\tau  < \infty.$$
Therefore we end up with the equality (2.16), where the kernel (2.17) can be written in the form

$$ S(x,\tau)  = {1\over \sqrt \pi} \int_0^\infty  e^{-y/2} h\left({x\over y} \right)  K_{i\tau} \left({y\over 2} \right)  {dy\over y}.\eqno(2.21)$$
Moreover, it can be calculated explicitly, substituting (2.20) in the right-hand side of (2.21) and employing relation  (3.35.9.3) in \cite{mel}. In fact, we get

$$ S(x,\tau)  =   { 3  \sqrt \pi \ \tau  \Gamma (-1- \mu)  \over  4 x  \  \Gamma(2- \mu) \sinh(\pi\tau) }  \ {}_3F_2 \left( 1+i\tau,\ 1-i\tau,\   {5\over 2} ;\  2+\mu,\ 2-\mu;\ - {1\over x} \right) $$

$$+ \left( {x\over 4}\right)^\mu  {   \sqrt \pi \ \Gamma (1+ \mu) \Gamma( -\mu -i\tau) \Gamma( -\mu +i\tau) \over   \Gamma(1-\mu) \Gamma( \mu-1/2) \Gamma(1/2-\mu) }$$

$$\times  \ {}_3F_2 \left( -\mu -i\tau,\ -\mu +i\tau, \ {3\over 2} - \mu ; \ - \mu,\ 1- 2\mu;\ - {1\over x} \right).\eqno(2.22)$$
Finally, returning to (2.16) and making a simple substitution, we come up with the following inversion formula for the index transform (1.1)

$$  f\left(x\right)  = - {1\over \pi^2}  \int_{1/x}^\infty \int_0^\infty \tau \sinh(2\pi\tau) S(y,\tau)  (F f)(\tau) {d\tau dy\over y}.\eqno(2.23)$$ 
One can change the order of integration in the right-hand side of (2.23) due to the estimate (see (2.19), (2.21)) 

$$\int_{1/x}^\infty \int_0^\infty \tau \sinh(2\pi\tau) \left| S(y,\tau)  (F f)(\tau) \right| {d\tau dy\over y}$$

$$\le C  ||Ff||_{L_1(\mathbb{R}_+;  \tau e^{3\pi \tau/2} d\tau)} \int_0^\infty e^{-u} u^{\nu-5/4} du  \int_{1/x}^\infty  y^{-\nu-1} dy < \infty, \ \nu > {1\over 4},$$
where $C >0$ is an absolute constant.   Hence from (2.22) and relation (1.16.1.1) in \cite{prud}, Vol. III

$$ \int_{1/x}^\infty  S(y,\tau) {dy\over y} = { 3  x \sqrt \pi \ \tau  \Gamma (-1- \mu)  \over  4  \  \Gamma(2- \mu) \sinh(\pi\tau) }  \ {}_4F_3 \left( 1+i\tau,\ 1-i\tau,\   {5\over 2},\ 1  ;\  2+\mu,\ 2-\mu,\ 2;\ -  x \right)$$

$$-  {1\over (4x)^\mu }   {   \sqrt \pi \ \Gamma (\mu) \Gamma( -\mu -i\tau) \Gamma( -\mu +i\tau) \over   \Gamma(1-\mu) \Gamma( \mu-1/2) \Gamma(1/2-\mu) }$$

$$\times  \ {}_3F_2 \left( -\mu -i\tau,\ -\mu +i\tau, \ {3\over 2} - \mu ; \ 1-\mu,\ 1- 2\mu;\ -  x \right),\eqno(2.24)$$
and after a simple substitution formula (2.23) can be written as follows

$$  f(x)  =  {1\over  \pi \sqrt \pi}   \int_0^\infty \left[ { \Gamma ( \mu)   \Gamma( -\mu -i\tau) \Gamma( -\mu +i\tau) \sinh(2\pi\tau) \over  (4 x)^\mu \   \Gamma(1-\mu) \Gamma( \mu-1/2) \Gamma(1/2-\mu)} \right.$$

$$\left. \times  \ {}_3F_2 \left( -\mu -i\tau,\ -\mu +i\tau, \ {3\over 2} - \mu ; \ 1-\mu,\ 1- 2\mu;\ -   x \right) \right.$$

$$\left. -  { 3 x \  \tau \cosh(\pi\tau)  \Gamma (-1- \mu)  \over  2 \  \Gamma(2- \mu)  }  \ {}_4F_3 \left( 1+i\tau,\ 1-i\tau,\   {5\over 2},\ 1  ;\  2+\mu,\ 2-\mu,\ 2;\ -  x  \right) \right] $$

$$\times (F f)(\tau)\  \tau d\tau .\eqno(2.25)$$ 
The hypergeometric functions in (2.25) can be expressed in terms of the associated Legendre functions. To do this, we,  firstly,  appeal to the differential property for the generalized hypergeometric function \cite{erd}, Vol. I to  have the formula

$$  {}_3F_2 \left( -\mu -i\tau,\ -\mu +i\tau, \ {3\over 2} - \mu ; \ 1-\mu,\ 1- 2\mu;\ -  x \right) = { x^{\mu+1/2}  \over 1/2-\mu}  $$

$$\times  {d\over dx} \left[  x^{1/2-\mu}  {}_3F_2 \left( -\mu -i\tau,\ -\mu +i\tau, \ {1\over 2} - \mu ; \ 1-\mu,\ 1- 2\mu;\ -  x \right) \right].$$
However,  the generalized hypergeometric function on the right-hand side of the latter equality can be written in terms of the squares of the associated Legendre functions (cf. \cite{prud}, Vol. III, Entries 7.4.1.26 and 7.3.1.40 ),  taking into account the fundamental identity  \cite{erd}, Vol. I

$$P^\mu_{\nu} (z)=  P^\mu_{-\nu-1} (z).\eqno(2.26)$$
Hence  we obtain

$${}_3F_2 \left( -\mu -i\tau,\ -\mu +i\tau, \ {1\over 2} - \mu ; \ 1-\mu,\ 1- 2\mu;\ -  x \right) = \left({x\over 4} \right)^{\mu}  {\Gamma^2(1-\mu) \over 2i\tau} $$

$$ \times  \left[ (i\tau +\mu)   \left[ P^{\mu}_{-i\tau} \left(\sqrt{1+x}\right) \right]^2  +   (i\tau -\mu)\    \left[ P^{\mu}_{i\tau} \left(\sqrt{1+x}\right) \right]^2 \right].\eqno(2.27)$$
In the same manner we write 

$${}_4F_3 \left( 1+i\tau,\ 1-i\tau,\   {5\over 2},\ 1  ;\  2+\mu,\ 2-\mu,\ 2;\ -  x \right) $$

$$= {2\over 3 \sqrt x}\  {d\over dx} \left[ x^{3/2} {}_4F_3 \left( 1+i\tau,\ 1-i\tau,\   {3\over 2},\ 1  ;\  2+\mu,\ 2-\mu,\ 2;\ -  x \right)\right].\eqno(2.28)$$
Moreover, in turn, (see \cite{prud}, Vol. III, Entry 1.16.1.1) 

$$ {}_4F_3 \left( 1+i\tau,\ 1-i\tau,\   {3\over 2},\ 1  ;\  2+\mu,\ 2-\mu,\ 2;\ -  x \right)$$

$$ =  (\mu-1 ) x^{\mu-1}  \int_x^\infty y^{- \mu}  {}_4F_3 \left( 1+i\tau,\ 1-i\tau,\   {3\over 2},\ 1  ;\  2+\mu,\ 1-\mu,\ 2;\ -  y\right) dy.\eqno(2.29)$$
But  the formula in \cite{prud}, Vol. III, Entry  7.5.1.5  says 

$$ {}_4F_3 \left( 1+i\tau,\ 1-i\tau,\   {3\over 2},\ 1  ;\  2+\mu,\ 1-\mu,\ 2;\ -  y\right) $$

$$=  {}_2F_1 \left( 1+i\tau,\ 1-i\tau;\   2+\mu;\ {1- \sqrt{1+y} \over 2}\right)  {}_2F_1 \left( 1+i\tau,\ 1-i\tau; \  1-\mu ;\   {1- \sqrt{1+y} \over 2}\right).$$
Besides,  via  the formula in \cite{prud}, Vol. III, Entry   7.3.1.92 we find the equalities 

$$ {}_2F_1 \left( 1+i\tau,\ 1-i\tau;\   1-\mu;\  {1- \sqrt{1+y} \over 2} \right) $$

$$=   {\Gamma (1-\mu)\over  i\tau}  \left({\sqrt{y}  \over \sqrt{1+y} +1 }\right)^{ \mu} 
\left[ (i\tau -\mu)  P^{\mu}_{-1-i\tau}(\sqrt{1+y}) -   {\sqrt{y}  \over \sqrt{1+y} +1 } P^{1+\mu}_{-1-i\tau}(\sqrt{1+y})\right],$$

$$ {}_2F_1 \left( 1+i\tau,\ 1-i\tau;\   2+\mu;\   {1- \sqrt{1+y} \over 2}\ \right) $$

$$=   {\Gamma (2+\mu)\over  i\tau}  \left({\sqrt{y}  \over \sqrt{1+y} +1 }\right)^{- (1+\mu)} 
\left[ (1+i\tau +\mu)  P^{-1-\mu}_{-1-i\tau}(\sqrt{1+y})\right.$$

$$\left.  -   {\sqrt{y}  \over \sqrt{1+y} +1 } P^{-\mu}_{-1-i\tau}(\sqrt{1+y})\right].$$
Finally, minding (2.26),  we get

$$ {}_4F_3 \left( 1+i\tau,\ 1-i\tau,\   {3\over 2},\ 1  ;\  2+\mu,\ 1-\mu,\ 2;\ -  y\right) =  -   {\Gamma (2+\mu) \Gamma (1-\mu)  \over  \tau^2 \sqrt y\  ( \sqrt{1+y} +1)}  $$

$$\times \left[ (1+i\tau +\mu) ( \sqrt{1+y} +1)   P^{-1-\mu}_{i\tau}(\sqrt{1+y})  -   \sqrt{y}\   P^{-\mu}_{i\tau}(\sqrt{1+y})\right]$$

$$\times \left[ (i\tau -\mu) ( \sqrt{1+y} +1)  P^{\mu}_{i\tau}(\sqrt{1+y}) -   \sqrt{y} \  P^{1+\mu}_{i\tau}(\sqrt{1+y})\right].$$

We summarize these results  as the following inversion   theorem. 

{\bf Theorem 2}. {\it Let ${\rm Re} \mu < 0,\ \mu \notin \mathbb{Z},\ \min(-  {\rm Re} \mu,\ 1) > 1/4$,  the  conditions of Theorem 1 hold and $\varphi \in  L_{3/4,1}(0,1) \cap L_{5/4,1}(1,\infty),\  f \in L_{-\nu,1}(\mathbb{R}_+), \  \min(-  {\rm Re} \mu,\ 1) > \nu >1/4,\   Ff  \in L_1(\mathbb{R}_+;   \tau  e^{3\pi \tau/2} d\tau)$ .  Let,  besides,  the Mellin transform $f^*(-\nu-it)$ be  such that   $(\nu+it) f^*(-\nu-it)  \in L_1( \mathbb{R})$.  Then for all $x >0$ the inversion formula $(2.25)$ holds. }

\section{Index transform (1.2)} 

In this section we investigate the boundedness properties of the  operator (1.2) and establish the inversion formula for this transformation.    We begin with 

{\bf Theorem 3.}  {\it Let ${\rm Re}\mu < 1/2$. The index transform  $(1.2)$  is well-defined as a  bounded operator $G: L_{1} \left(\mathbb{R}_+\right) \to L_{\nu,\infty}  (\mathbb{R}_+),\  \nu \in ( {\rm Re}\mu,\ 1/2)$ and the following norm inequality takes place
$$||G g||_{\nu,\infty}  \equiv \hbox{ess sup}_{x >0} | x^\nu (G g)(x)| \le C_{\mu,\nu}  ||g ||_{L_1(\mathbb{R}_+)},\eqno(3.1)$$
where $C_{\mu,\nu}$ is defined by $(2.2)$.   Moreover,   if, besides,  $G g \in L_{1-\gamma,1}(\mathbb{R}_+),\ 1/2 < \gamma < 1- {\rm Re}\mu$, $G^* g$ is analytic in the strip ${\rm Re}\mu < {\rm Re} s < 1- {\rm Re}\mu$ and 

$$\sup_{\mu_0 \le  {\rm Re} s \le \mu_1}  \left| (G^* g) \left(1- {\rm Re} s+ it\right) \right| \to 0,\quad |t| \to \infty,\ [\mu_0, \mu_1] \subset ({\rm Re}\mu ,  1- {\rm Re}\mu),\eqno(3.2)$$
where $\max\left( {\rm Re}\mu, 0\right)   < \mu_0 < 1/2 < \mu_1$, then for all $y >0$ }
$${1\over 2\pi i}   \int_{\nu -i\infty}^{\nu  +i\infty} \frac{\Gamma(s) \Gamma (s-\mu) \Gamma(1/2-s)} { \Gamma(s-1/2)  \Gamma(1-s-\mu) } (G^* g) (1-s) y^{ -s} ds$$

$$=  \sqrt {\pi} \ e^{y/2}   \int_{0}^\infty   \   K_{i\tau} \left({y\over 2} \right)  {g(\tau) \over \cosh(\pi\tau)} d\tau, \quad \nu \in [ \mu_0, 1/2 ).\eqno(3.3)$$

\begin{proof}   Indeed, as we see from (1.4) the kernel $ \Phi(x,\tau )$  has a bound

$$|\Phi(x,\tau )| \le C_{\mu,\nu}\  x^{-\nu},\ x >0,$$
where $C_{\mu,\nu}$ is defined by (2.2), we have 

$$|(G g)(x)| \le C_{\mu,\nu} x^{-\nu} \int_{\mathbb{R}_+} |g(\tau)| d\tau =  C_{\mu,\nu} x^{-\nu} ||g ||_{L_1(\mathbb{R}_+)},$$
and (3.1) holds.   Further, taking the Mellin transform (1.6) of both sides of (1.2) at the   point $1-s$ under the condition  $Gg \in L_{1-\gamma,1}(\mathbb{R}_+), \ 1/2 < \gamma < 1- {\rm Re}\mu$, we change the order of integration on the right-hand side of the obtained equality by Fubini's theorem. Moreover, taking into account (1.4), we derive 

$$  \frac{\Gamma(s) \Gamma (s-\mu) } { \Gamma(s-1/2)\Gamma(1-s-\mu) } (G^* g) (1-s)=   \int_{0}^\infty \Gamma\left(s + i\tau\right) \Gamma\left(s - i\tau \right)  g(\tau) d\tau.\eqno(3.4)$$
Meanwhile,  employing relation (8.4.23.5) in \cite{prud}, Vol. III

$${ \sqrt \pi\over \cosh(\pi\tau)} \ e^{x/2}  K_{i\tau} \left({x\over 2}\right)  = {1\over 2\pi i}   \int_{\nu -i\infty}^{\nu  +i\infty}   {\Gamma\left( s + i\tau\right) \Gamma\left(s - i\tau \right) \Gamma( 1/2- s) } x^{-s} ds,\eqno(3.5)$$
we multiply  (3.4)  by $\Gamma(1/2-s)$ and then take  the inverse Mellin transform (1.8) of both sides of the obtained equality over the vertical line $(\nu-i\infty, \nu +i\infty),   \nu \in [ \mu_0, 1/2 )$. This is indeed possible via (3.2) and conditions of the theorem. Thus we establish (3.3), where the integral in the left-hand side converges absolutely.  Theorem 3 is proved. 
\end{proof}

The inversion formula for the index transform (1.2) is given by

{\bf Theorem 4}.  {\it Let ${\rm Re}\mu < 1/2,\ \mu \notin \mathbb{Z}$,  $g(z/i)$ be an even analytic function in the strip $D= \left\{ z \in \mathbb{C}: \ |{\rm Re} z | \le  \alpha < 1/2\right\}$,  such  that $ g(0)=g^\prime (0)=0$ and $g(z/i)$ be  absolutely  integrable over any vertical line in  $D$.   If $Gg \in  L_{1- \nu, 1} \left(\mathbb{R}_+\right),\   \nu \in [ \mu_0, 1/2 ),  \mu_0 \in ( \max\left( {\rm Re}\mu, 0\right) , 1/2) $, then for all  $x \in \mathbb{R}_+$ the  inversion formula holds for the index transform $(1.2)$} 

$$  g(x)  = \lim_{\varepsilon \to 0+}  { 1\over \pi^2} \ x\sinh (2\pi x)\   \int_0^\infty \left[    \frac {\sqrt \pi\  \Gamma(\varepsilon -ix) \Gamma(\varepsilon +ix) }{ 2\mu\ \Gamma(1/2+\varepsilon) } \right.$$

$$\left. \times  {}_4F_3 \left( {1\over 2},\ {3\over 2},\ \varepsilon -ix,\ \varepsilon + ix;\ 1+\mu,\ 1-\mu,\ {1\over 2}+\varepsilon ; \ -u\right)   \right.$$

$$\left. +   \frac {\pi\  4^\mu  \Gamma(\mu)   \Gamma(\varepsilon -\mu-ix) \Gamma(\varepsilon-\mu  +ix) }{ \Gamma(1/2+\varepsilon-\mu)  \Gamma(\mu-1/2) \Gamma(1-\mu)} \right.$$

$$\left. \times \ {}_4F_3 \left( {1\over 2} - \mu,\ {3\over 2} - \mu, \ \varepsilon-\mu -ix,\ \varepsilon-\mu + ix;\ 1 - \mu,\ 1- 2\mu, \ {1\over 2}+\varepsilon-\mu ; \ -u\right) \right]  $$

$$\times (Gf)(u)  du.\eqno(3.6)$$

\begin{proof}    Indeed,  recalling  (3.3),  we write its left-hand side, appealing to the Parseval equality (1.7), to get

$${1\over 2\pi i}   \int_{\nu -i\infty}^{\nu  +i\infty} \frac{\Gamma(s) \Gamma (s-\mu) \Gamma(1/2-s)} { \Gamma(s-1/2)  \Gamma(1-s-\mu) } (G^* g) (1-s) y^{ -s} ds$$

$$ = \int_0^\infty U_\mu(y u) (Gf)(u) du,\eqno(3.7)$$
where

$$U_\mu(y) = {1\over 2\pi i}  \int_{\nu -i\infty}^{\nu  +i\infty} \frac{\Gamma(s) \Gamma (s-\mu) \Gamma(1/2-s)} { \Gamma(s-1/2)  \Gamma(1-s-\mu) } y^{-s} ds,\ y >0.$$
The function $U_\mu(y),\ \mu \notin \mathbb{Z}$ can be calculated explicitly by Slater's theorem in terms of a combination of two hypergeometric functions ${}_2F_2$. Precisely, we obtain

$$U_\mu(y) =  {1\over 2 \mu} \  {}_2F_2 \left( {1\over 2},\ {3\over 2};\ 1+\mu,\ 1-\mu; \ -y\right) $$

$$+ \left( {y\over 4}\right)^{-\mu} \frac {\sqrt \pi\  \Gamma(\mu) }{ \Gamma(\mu-1/2) \Gamma(1-\mu)} \ {}_2F_2 \left( {1\over 2} - \mu,\ {3\over 2} - \mu;\ 1 - \mu,\ 1- 2\mu; \ -y\right).$$ 
Returning to (3.3), we substitute its left-hand side by the right-hand side of (3.7) and multiply  both sides of the obtained equality by $ e^{-y/2} K_{ix} \left({y/2} \right) y^{\varepsilon - 1}$ for some positive $\varepsilon \in (0,1)$. Then,  integrating  with respect to $y$ over $(0, \infty)$, we appeal to the asymptotic behavior of the hypergeometric functions at infinity \cite{erd}, Vol. I and the condition $Gg \in  L_{1- \nu, 1} \left(\mathbb{R}_+\right), \nu \in [ \mu_0, 1/2 )$ to justify  the interchange of the order of integration in the left-hand side  due to the absolute convergence of the  iterated integral.   Finally, the inner integral is calculated  with the use of formula 3.35.9.3 in \cite{mel}, and we derive  the equality 

$$\int_0^\infty \int_0^\infty  e^{-y/2} K_{ix} \left({y/2} \right) y^{\varepsilon - 1} U_\mu(y u)  (Gf)(u) dy du $$

$$=   \int_0^\infty \left[    \frac {\sqrt \pi\  \Gamma(\varepsilon -ix) \Gamma(\varepsilon +ix) }{ 2\mu\ \Gamma(1/2+\varepsilon) }  {}_4F_3 \left( {1\over 2},\ {3\over 2},\ \varepsilon -ix,\ \varepsilon + ix;\ 1+\mu,\ 1-\mu,\ {1\over 2}+\varepsilon ; \ -u\right)   \right.$$

$$\left. +   \frac {\pi\  4^\mu  \Gamma(\mu)   \Gamma(\varepsilon -\mu-ix) \Gamma(\varepsilon-\mu  +ix) }{ \Gamma(1/2+\varepsilon-\mu)  \Gamma(\mu-1/2) \Gamma(1-\mu)} \right.$$

$$\left. \times \ {}_4F_3 \left( {1\over 2} - \mu,\ {3\over 2} - \mu, \ \varepsilon-\mu -ix,\ \varepsilon-\mu + ix;\ 1 - \mu,\ 1- 2\mu, \ {1\over 2}+\varepsilon-\mu ; \ -u\right) \right]  $$

$$\times (Gf)(u)  du =  \int_0^\infty   K_{ix} \left({y\over 2} \right) y^{\varepsilon - 1}  \int_{0}^\infty   \   K_{i\tau} \left({y\over 2} \right)  {g(\tau) \over \cosh(\pi\tau)} d\tau dy.\eqno(3.8)$$
In the meantime,  the right-hand side of latter equality in (3.8) can be treated, using   the evenness of $g$ and  the representation  of the modified Bessel  function $K_z(y)$ in terms of the modified Bessel function of the first kind $I_z(y)$ \cite{erd}, Vol. II. Hence with a simple substitution we find 
$$ {1\over 2}  \int_0^\infty  K_{ix} \left({y\over 2} \right) y^{\varepsilon -1} \int_{-\infty}^\infty  K_{i\tau} \left({y\over 2} \right)  {g(\tau) \over \cosh(\pi\tau) } d\tau dy$$

$$=  \pi i \int_0^\infty  K_{ix} \left({y\over 2} \right) y^{\varepsilon -1} \int_{-i\infty}^{i\infty}   I_{ z} \left({y\over 2} \right)  {g(z/i) \over \sin (2\pi z) } dz\  dy. \eqno(3.9)$$
On the other hand, according to our assumption $g(z/i)$ is analytic and integrable in the vertical  strip $0\le  {\rm Re}  z \le \alpha< 1/2$,  $g(0)=g^\prime (0)=0$.  Hence,  appealing to the inequality for the modified Bessel   function of the first  kind  (see \cite{yal}, p. 93)
 $$|I_z(y)| \le I_{  {\rm Re} z} (y) \  e^{\pi |{\rm Im} z|/2},\   0< {\rm Re} z \le \alpha,$$
one can move the contour to the right in the latter integral in (3.8). Then 

$$ \pi i  \int_0^\infty  K_{ix} \left({y\over 2} \right) y^{\varepsilon -1} \int_{-i\infty}^{i\infty}   I_{ z} \left({y\over 2} \right)  {g(z/i) \over \sin (2\pi z) } dz\  dy$$

$$=  \pi i    \int_0^\infty  K_{ix} \left({y\over 2} \right) y^{\varepsilon -1} \int_{\alpha -i\infty}^{\alpha + i\infty}   I_{ z} \left({y\over 2} \right)  {g(z/i) \over \sin (2\pi z) } dz\  dy.$$
Now ${\rm Re} z >0$,  and  it is possible to pass to the limit under the integral sign when $\varepsilon \to 0$ and to change the order of integration due to the absolute and uniform convergence.  Therefore the value of the integral (see relation (2.16.28.3) in \cite{prud}, Vol. II)
$$\int_0^\infty K_{ix}(y) I_z(y) {dy\over y} = {1\over x^2 + z^2} $$ 
leads us to the equalities 

$$\lim_{\varepsilon \to 0}   \pi i  \int_0^\infty  K_{ix} \left({y\over 2} \right) y^{\varepsilon -1} \int_{-i\infty}^{i\infty}   I_{ z} \left({y\over 2} \right)  {g(z/i) \over \sin (2\pi z) } dz\  dy$$

$$=     \pi i  \int_{\alpha -i\infty}^{\alpha + i\infty}   {g(z/i) \over (x^2+ z^2) \sin (2\pi z) } dz $$

$$=  {\pi   i\over 2}    \left( \int_{-\alpha +i\infty}^{- \alpha- i\infty}   +   \int_{\alpha -i\infty}^{ \alpha+  i\infty}   \right)  {  g(z/i) \  dz \over (z-ix) \  z \sin(2\pi z)}. \eqno(3.10)$$
Hence conditions of the theorem allow to apply the Cauchy formula in the right-hand side of the latter equality in (3.8).  Thus 
$$\lim_{\varepsilon \to 0}   \pi i   \ \int_0^\infty  K_{ix} \left({y\over 2} \right) y^{\varepsilon -1} \int_{-i\infty}^{i\infty}   I_{ z} \left({y\over 2} \right)  {g(z/i) \over \sin (2\pi z) } dz\  dy$$

$$ =  { \pi^{2} \  g(x) \over  x\sinh (2\pi x)} ,\quad x > 0.\eqno(3.11)$$
Therefore passing to the limit through  equalities  (3.8), we end up with the inversion formula (3.6), completing the proof of Theorem 4. 

 \end{proof}
 
 {\bf Remark 1}.  {\it If the passage to the limit in $(3.6)$ is permitted, the inversion formula takes the form}
 
 $$  g(x)  =  { 1\over \pi\mu}    \int_0^\infty \left[   \cosh(\pi x)   {}_3F_2 \left(  {3\over 2},\  -ix,\   ix;\ 1+\mu,\ 1-\mu ; \ -u\right)   \right.$$

$$\left. +   \frac {  4^\mu  \Gamma(1+\mu)   \Gamma( -\mu-ix) \Gamma(-\mu  +ix) x\sinh (2\pi x)}{ \Gamma(1/2-\mu)  \Gamma(\mu-1/2) \Gamma(1-\mu)} \right.$$

$$\left. \times \ {}_3F_2 \left(  {3\over 2} - \mu, \ -\mu -ix,\ -\mu + ix;\ 1 - \mu,\ 1- 2\mu ; \ -u\right) \right]   (Gf)(u)  du.$$

\section{Boundary   value problem}

In this section  the index transform (1.2) is employed to investigate  the  solvability  of the boundary  value  problem  for the  following third   order partial differential  equation

$$ 2x \left( y^2+ x^2 \left(1+{1\over r}\right) \right)  {\partial^3 u \over \partial x^3} +   2y \left( x^2+ y^2 \left(1+{1\over r}\right) \right)  {\partial^3 u \over \partial y^3}$$

$$ +  2x \left( x^2+ y^2 \left(1+{3\over r}\right) \right)  \ {\partial^3 u \over \partial x\partial y^2} +  2y \left( y^2+ x^2 \left(1+{3\over r}\right) \right)   {\partial^3 u \over \partial y \partial x^2}$$

$$+   \left( 5 y^2+ 3 x^2 \left(3+{2\over r}\right) \right)  \  {\partial^2 u \over \partial x^2} +   \left( 5 x^2+ 3 y^2 \left(3+{2\over r}\right) \right)\   {\partial^2 u \over \partial y^2}$$

$$  + 4xy \left( 2+{3\over r}\right) \  {\partial^2 u \over \partial x\partial y}  + 2  \left( 4+  {1-\mu^2 \over r}  \right)  \left( x {\partial u \over \partial x}  + y {\partial u  \over \partial y}\right)  + u  =0,\eqno(4.1)$$
where $(x,y) \in \mathbb{R}^2 \backslash \{ 0\},\  r= \sqrt {x^2+y^2}$.    Writing  (4.1) in polar coordinates $(r,\theta)$, we end up with the equation

$$2 r^2(1+r)  {\partial^3  u \over \partial r^3} +  2r  {\partial^3  u \over \partial r \partial \theta^2} +  r \left( 11 r+  6\right) {\partial^2 u  \over \partial r^2} $$

$$ +  {\partial^2  u \over  \partial \theta^2} +   \left( 2(1-\mu^2) +   11 r \right)  {\partial u \over \partial r} + u = 0.\eqno(4.2)$$

{\bf Lemma 3.} {\it  Let $0 < {\rm Re}\  \mu < 1/2,\  g  \in L_1\left(\mathbb{R}_+\right),\  \beta \in (0, 2\pi)$. Then  the function
$$u(r,\theta)=   \sqrt{ {\pi\over 1+r}}  \int_0^\infty \Gamma(1+i\tau-\mu)  \Gamma(1-i\tau-\mu)  $$

$$\times  P^\mu_{i\tau} (\sqrt{1+r}) P^\mu_{- i\tau} (\sqrt{1+r})  {\sinh(\theta \tau)\over \sinh(\beta\tau)} g(\tau) d\tau\eqno(4.3)$$
 satisfies   the partial  differential  equation $(4.2)$ on the wedge  $(r,\theta): r   >0, \  0\le \theta < \beta$, vanishing at infinity.}

\begin{proof} The proof  is straightforward by  substitution (4.3) into (4.2) and the use of  (1.15).  The necessary  differentiation  with respect to $r$ and $\theta$ under the integral sign is allowed via the absolute and uniform convergence, which can be verified, appealing to  the integrability condition $g \in L_1\left(\mathbb{R}_+\right),\  \beta \in (0, 2\pi)$ and estimates of  derivatives of the kernel (1.3) with respect to $r$.  Finally,  the condition $ u(r,\theta) \to 0,\ r \to \infty$  is due to the estimate of the integral (1.4).  
\end{proof}

Finally  we will formulate the boundary  value problem for equation (4.2) and give its solution.

{\bf Theorem 5.} {\it Let  $g(x)$ be given by formula $(3.6)$ and its transform $(Gg) (y)\equiv G(y)$ satisfies conditions of Theorem $4$.  Then  $u (r,\theta),\   r >0,  \  0\le \theta < \beta$ by formula $(4.3)$  will be a solution  to the boundary  value problem on the wedge for the partial differential  equation $(4.2)$ subject to  boundary  conditions}
$$u(r,0) = 0,\quad\quad  u(r,\beta) = G(r).$$

\bigskip
\noindent{{\bf Disclosure statement}}
\bigskip

\noindent  No potential conflict of interest was reported by the author. 

\bigskip
\noindent{{\bf Funding}}
\bigskip

\noindent The work was partially supported by CMUP (UID/MAT/00144/2019),  which is funded by FCT (Portugal) with national (MEC),   European structural funds through the programs FEDER  under the partnership agreement PT2020, and Project STRIDE - NORTE-01-0145-FEDER- 000033, funded by ERDF - NORTE 2020.

\bibliographystyle{amsplain}

\begin{thebibliography}{10}


\bibitem{yak}   Yakubovich S.  Index transforms.   Singapore:  World Scientific Publishing Company; 1996.

\bibitem{erd}    Erd\'elyi A,  Magnus W,   Oberhettinger  F,   Tricomi FG.  Higher transcendental functions. Vols. I,  II. New  York: McGraw-Hill;  1953.

\bibitem{vir}   Virchenko N, Fedotova I.  Generalized associated Legendre functions and their applications (with a Foreword by Semyon Yakubovich).   Singapore:  World Scientific Publishing Company;  2001.

\bibitem{prud}  Prudnikov AP,  Brychkov  YuA,  Marichev OI. Integrals and series:  Vol. I: Elementary functions. New York:  Gordon and Breach;   1986;   Vol. II:  Special functions. New York: Gordon and Breach;  1986;   Vol. III:  More special functions. New York:   Gordon and Breach; 1990.

\bibitem{yal}   Yakubovich S,  Luchko Yu.  The hypergeometric approach to integral transforms and convolutions, Mathematics and its applications.  Vol. 287.  Dordrecht:  Kluwer Academic Publishers Group; 1994.

\bibitem {tit}  Titchmarsh EC.   An introduction to the theory of Fourier integrals.   New York:  Chelsea; 1986.


\bibitem{square}  Lebedev NN.   On an integral representation of an arbitrary function in terms of squares of Macdonald functions with imaginary index,  {\it Sibirsk. Mat. Zh.},   {\bf 3}  (1962),  213-222 (in Russian).


\bibitem{mel}   Brychkov  YuA,  Marichev OI, Savischenko NV. Handbook of Mellin transforms.   Boca Raton: CRC Press; 2018. 











\end{thebibliography}

\end{document}